# RISK HULL METHOD AND REGULARIZATION BY PROJECTIONS OF ILL-POSED INVERSE PROBLEMS

By L. Cavalier and Yu. Golubev

*Université de Provence (Aix-Marseille 1)*

We study a standard method of regularization by projections of the linear inverse problem $Y = Af + \epsilon$, where $\epsilon$ is a white Gaussian noise, and $A$ is a known compact operator with singular values converging to zero with polynomial decay. The unknown function $f$ is recovered by a projection method using the singular value decomposition of $A$. The bandwidth choice of this projection regularization is governed by a data-driven procedure which is based on the principle of risk hull minimization. We provide nonasymptotic upper bounds for the mean square risk of this method and we show, in particular, that in numerical simulations this approach may substantially improve the classical method of unbiased risk estimation.

**1. Introduction and main result.** The inverse problem paradigm is related to the classical linear algebra problem in which we want to find a solution $x \in \mathbb{R}^d$ of the linear equation

$$Ax = y, \tag{1.1}$$

where $A$ is a known $d \times d$ matrix and $y$ is a given vector in $\mathbb{R}^d$. From a mathematical viewpoint, the linear inverse problem can be considered a straightforward generalization of (1.1). Let $\mathbb{H}$, $\mathbb{G}$ be two Hilbert spaces and let $A$ be a continuous linear operator $\mathbb{H} \to \mathbb{G}$. Suppose we have at our disposal an element (a function) defined by

$$Y = Af + \epsilon, \tag{1.2}$$

where $\epsilon$ is an unknown function which is small. The goal is to recover $f \in \mathbb{H}$.

Numerous applications of inverse problems in medical image processing, econometrics and astrophysics make this area very attractive for mathematical incursions. The mathematical literature on inverse problems is so vast that it would be impractical to cite it here. We refer the interested readers to [3, 10, 14], where interesting applications of inverse problems can be found.

In the last two decades, the stochastic approach, which goes back to [18], has been very intensively studied in the statistical literature (see, e.g., [5, 7, 8, 9, 11, 16, 17, 19, 21]). In this approach, it is usually assumed that $\epsilon$ is a Gaussian white noise in $\mathbb{H}$ (see, for details, [15]).

The simplest way to understand why the problem (1.2) may be difficult is to look at the singular value decomposition (SVD) of $A$. Let $A^*$ be the adjoint to $A$. Suppose $A^*A$ is a compact operator with eigenvalues $\lambda_k \geq 0, k = 1, \ldots,$ and eigenfunctions $\varphi_k, k = 1, \ldots.$ Let $\psi_k = A\varphi_k/\|A\varphi_k\|$. Then we get the following equivalent representation of (1.2):

$$(1.3) \qquad y_k = \theta_k + \sigma_k \xi_k, \qquad k = 1, 2, \ldots,$$

where $\xi_k$ are i.i.d. $\mathcal{N}(0,1)$, $y_k = \langle Y, \psi_k \rangle / \sqrt{\lambda_k}, \theta_k = \langle f, \varphi_k \rangle, \sigma_k = \varepsilon/\sqrt{\lambda_k}$, and $\varepsilon$ is a known spectral density of Gaussian white noise $\epsilon$.

Ill-posed inverse problems are characterized by the fundamental property that $\sigma_k \to \infty$ as $k \to \infty$, and the behavior of $\sigma_k$ for large $k$ describes the difficulty of the inverse problem. In this paper we will deal with moderately ill-posed inverse problems with polynomially increasing ($\sigma_k \asymp k^\beta, \beta \geq 0$). Recall that in the statistical literature this type of inverse problem is often associated with estimation of the derivative of order $\beta$ of a regression function.

The fact that $\sigma_k \to \infty$ immediately entails that the natural inversion

$$(1.4) \qquad A^{-1}Y = \sum_{k \,:\, \lambda_k > 0} \lambda_k^{-1/2} \langle Y, \psi_k \rangle \varphi_k$$

cannot be used since the quadratic risk of this method is infinite. A standard way to overcome this difficulty is based on a regularization technique. Nowadays the family of regularization methods available for practical applications is very large; see [10] and [23]. In the present paper, we will focus on regularization by projections. The idea of this method is very simple. In order to invert $A$, let us use the first $N$ terms of the expansion (1.4). In other words, to recover $f$ or equivalently $\theta_k, k = 1, \ldots,$ in the model (1.3), we use the projection method

$$(1.5) \qquad \tilde{\theta}_k(N) = y_k \mathbf{1}(k \leq N).$$

The mean square risk of this inversion method is computed very easily:

$$(1.6) \qquad R(\theta, N) = \mathbf{E}_\theta \|\tilde{\theta}(N) - \theta\|^2 = \sum_{k=N+1}^{\infty} \theta_k^2 + \sum_{k=1}^{N} \sigma_k^2.$$

The parameter $N$ here is called the bandwidth and the major statistical problem is related to the data-driven choice of $N$. Roughly speaking, the goal of this choice is to minimize the right-hand side of (1.6) based on the noisy data $y_k$ from (1.3).

A classical approach to this minimization problem is based on the principle of *unbiased risk estimation* (URE) (see [22]). The idea to use this method for adaptive bandwidth choice goes back to [1] and [20]. Originally, URE was proposed in the context of regression estimation $\sigma_k = \varepsilon$. Nowadays, it is used as a basic adaptation tool for many statistical models. For inverse problems, this method was studied in [5], where precise oracle inequalities for the mean square risk were obtained.

The heuristic motivation of URE is rather simple. The underlying optimization problem can be reformulated as minimization of $-\sum_{k=1}^{N} \theta_k^2 + \sum_{k=1}^{N} \sigma_k^2$ [see (1.6)]. Noticing that the unobservable term $\sum_{k=1}^{N} \theta_k^2$ can be estimated by $\sum_{k=1}^{N} (y_k^2 - \sigma_k^2)$, we choose the bandwidth as

$$(1.7) \quad N_{\mathrm{ure}}(y) = \arg\min_{N \geq 1} \bar{R}(y, N) \qquad \text{where } \bar{R}(y, N) = \left\{ -\sum_{k=1}^{N} y_k^2 + 2\sum_{k=1}^{N} \sigma_k^2 \right\}.$$

Intuitively, since $N_{\mathrm{ure}}(y)$ minimizes the estimator of the risk, it means that the risk of the method $\mathbf{E}_\theta \|\tilde{\theta}(N_{\mathrm{ure}}) - \theta\|^2$ can be controlled by the risk of the best projection method $\inf_N R(\theta, N)$, which is sometimes called *risk of oracle*. Following [4], we measure the quality of the method $\tilde{\theta}(N_{\mathrm{ure}})$ by the ratio of its risk to the risk of oracle,

$$(1.8) \quad r(\theta) = \frac{\mathbf{E}_\theta \|\tilde{\theta}(N_{\mathrm{ure}}) - \theta\|^2}{\inf_N R(\theta, N)}.$$

When we use URE we hope that $r(\theta)$ is bounded from above by a relatively small constant uniformly over all $\theta$. It is well known that this native hypothesis holds (see [4]) for direct estimation ($\sigma_k \equiv \varepsilon$). However, when we deal with an inverse problem the situation becomes more difficult.

In order to illustrate the difference between direct and inverse estimation, we will carry out a very simple numerical experiment. Obviously, we cannot compute in a numerical experiment $r(\theta)$ for all $\theta \in \mathbf{l}_2$. Therefore, let us take $\theta_k \equiv 0$ and compute $r(0)$ for two cases, $\sigma_k \equiv \varepsilon$ and $\sigma_k = \varepsilon k$. The first case corresponds to classical regression function estimation (direct estimation), whereas the second is related to the estimation of the first-order derivative of a regression function. Notice that in both cases the risk of the oracle is evidently $\inf_N R(0, N) = \varepsilon^2$ since $\arg\min_N R(0, N) = 1$. In order to shed some light on the performance of URE, we generated 2000 independent random vectors $y^j, j = 1, \ldots, 2000$, with the components defined by (1.3). For each vector we computed $N_{\mathrm{ure}}(y^j)$ and the normalized error $\|\hat{\theta}[N_{\mathrm{ure}}(y^j)] - \theta\|^2/\varepsilon^2$ and plotted these values as a stem diagram. We also computed the

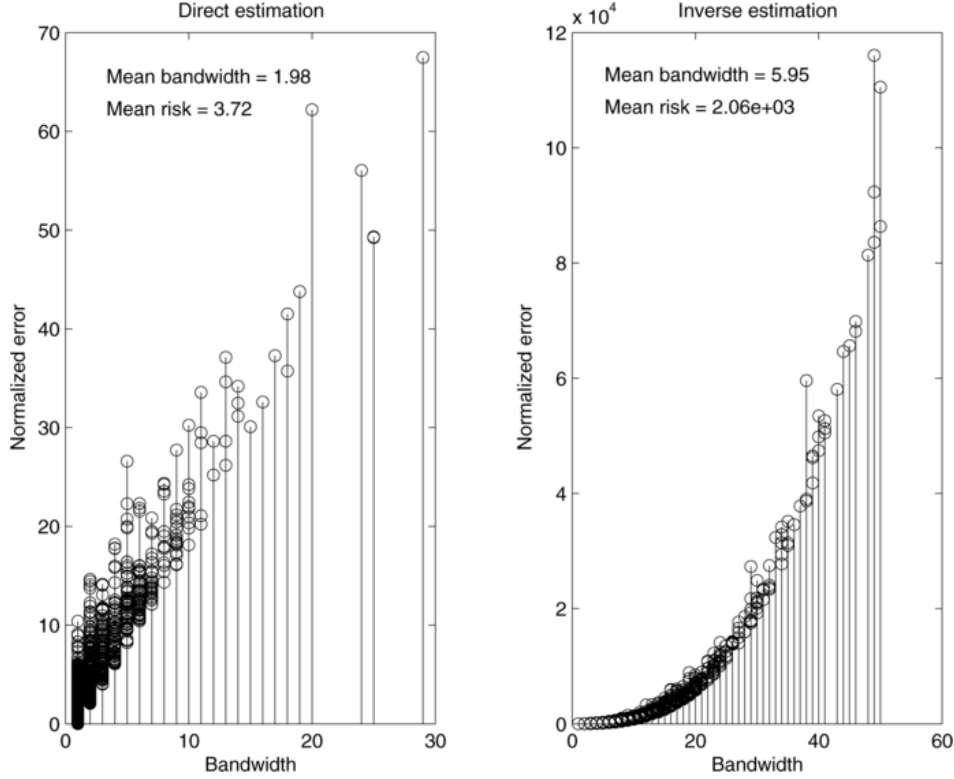

Fig. 1. *The method of unbiased risk estimation.*

mean empirical bandwidth $N_{\text{emp}}$ and the normalized mean empirical risk $R_{\text{emp}}$ by

$$N_{\text{emp}} = \frac{1}{2000}\sum_{j=1}^{2000} N_{\text{ure}}(y^j), \qquad R_{\text{emp}} = \frac{1}{2000\varepsilon^2}\sum_{j=1}^{2000} \|\hat{\theta}[N_{\text{ure}}(y^j)] - \theta\|^2.$$

Let us discuss briefly the numerical results of this experiment shown in Figure 1. The first display (direct estimation) shows that the URE method works reasonably well. Almost all bandwidths $N_{\text{ure}}(y^j)$ are relatively small (their mean is 1.98) and $r(0) = 3.72$. Even a quick look at the second display shows that the distribution of $N_{\text{ure}}(y^j)$ changed essentially. Now the mean is 5.95 and there are sufficiently many bandwidths $N_{\text{ure}}(y^j)$ greater than 20. This results in a catastrophic $r(0) \approx 2000$. On the other hand, it follows from the oracle inequalities (see [5, 6, 7] or Theorem 4 of the present paper) that in both cases there exist a lot of $\theta$ for which $r(\theta) \approx 1$. Comparing this fact with the simulations, we can conclude that for ill-posed inverse problems, URE does not work properly since very large $r(0)$ undermines its basic idea.

There exists a more general approach which is very close to URE. This method is called *method of penalized empirical risk*, and in the context of our problem it provides us with the bandwidth choice

$$N(y) = \arg\min_{N \geq 1} \bar{R}_{\text{pen}}(y, N),$$

(1.9)

$$\bar{R}_{\text{pen}}(y, N) = \left\{ -\sum_{k=1}^{N} y_k^2 + \sum_{k=1}^{N} \sigma_k^2 + \text{pen}(N) \right\},$$

where $\text{pen}(N)$ is a penalty function. The modern literature on this method is vast and we refer the interested reader to [2] or [4]. The main idea at the heart of this approach is that severe penalties permit one to improve substantially the performance of URE. For instance, it is known that this approach works well for severely ill-posed problems, where URE completely fails (see, e.g., [12]). However, it should be mentioned that the principal difficulty of this method is related to the choice of the penalty function $\text{pen}(N)$.

In this paper we propose a more general approach, called *risk hull minimization* (RHM), which gives a relatively good strategy for the penalty choice. Our goal is to present heuristic and mathematical justifications of this method. In the framework of the empirical risk minimization RHM can be defined as follows. Let the penalty in (1.9) be

(1.10) $$\text{pen}(N) = \text{pen}_{\text{rhm}}(N) = \sum_{k=1}^{N} \sigma_k^2 + (1 + \alpha)U_0(N),$$

where

$$U_0(N) = \inf\{t > 0 : \mathbf{E}\eta_N I(\eta_N \geq t) \leq \sigma_1^2\} \qquad \text{with } \eta_N = \sum_{i=1}^{N} \sigma_i^2(\xi_i^2 - 1).$$

(1.11)

RHM chooses the bandwidth $N_{\text{rhm}}(y) = N(y)$ according to (1.9) with the penalty function defined by (1.10), (1.11). The following theorem provides an upper bound for the mean square risk of this approach. Recall that it assumed that $\sigma_k$ has polynomial growth ($\sigma_k \asymp \varepsilon k^\beta$); see, for details, (2.5) and (2.6).

THEOREM 1. *There exist constants $C_* > 0$ and $\gamma_0 > 0$ such that for all $\gamma \in (0, \gamma_0]$ and $\alpha > 1$*

(1.12) $$\mathbf{E}\|\tilde{\theta}(N_{\text{rhm}}) - \theta\|^2 \leq (1 + \gamma)\inf_N R_{\text{rhm}}(\theta, N) + C_* \sigma_1^2 \left( \frac{1}{\gamma^{4\beta+1}} + \frac{1}{\alpha - 1} \right),$$

where $R_{\text{rhm}}(\theta, N) = \sum_{k=N+1}^{\infty} \theta_k^2 + \sum_{k=1}^{N} \sigma_k^2 + (1 + \alpha)U_0(N)$.

The statistical sense of Theorem 1 is rather transparent. The principal term of this upper bound is $\inf_N R_\alpha(\theta, N)$. The residual term

$$\min_\gamma \left\{ \gamma \inf_N R_{\mathrm{rhm}}(\theta, N) + \frac{C^* \sigma_1^2}{\gamma^{4\beta+1}} \right\} + \frac{C^* \sigma_1^2}{\alpha - 1}$$
$$= C^* \sigma_1^2 [(4\beta + 1)^{1/(4\beta+2)} + (4\beta + 1)^{-(4\beta+1)/(4\beta+2)}]$$
$$\times \left[ \frac{\inf_N R_{\mathrm{rhm}}(\theta, N)}{C^* \sigma_1^2} \right]^{(4\beta+1)/(4\beta+2)} + \frac{C^* \sigma_1^2}{\alpha - 1}$$

defines how much we should pay for stochastic minimization. Using this theorem we can get a typical panorama of minimax facts related to moderately ill-posed problems (see [5]). Moreover, simulations in Section 3 reveal that the constant $C_*$ is really small. It means, in particular, that in contrast to URE this method is stable.

The present paper is organized as follows. In Section 2, a heuristic motivation and additional facts related to RHM are presented. Section 3 contains simulation results. The proofs and technical lemmas are postponed to Section 4.

## 2. The RHM method.

2.1. *A heuristic motivation.* The heuristic motivation of the RHM approach is based on the oracle ideology. Suppose there is an oracle which provides us with $\theta_k, k = 1, \ldots,$ but we are allowed to use only projection methods. In this case the optimal bandwidth is evidently given by

$$N_{\mathrm{or}} = \arg\min_N r(y, N) \qquad \text{where } r(y, N) = \|\tilde{\theta}(N) - \theta\|^2.$$

Let us try to mimic this bandwidth choice. At the first glance this problem seems hopeless since in the decomposition

$$r(y, N) = \sum_{k=N+1}^{\infty} \theta_k^2 + \sum_{k=1}^{N} \sigma_k^2 \xi_k^2,$$

neither $\theta_k^2$ nor $\xi_k^2$ is really known. However, suppose for a moment that we know all the $\theta_k^2$, and we try to minimize $r(y, N)$. Since $\xi_k^2$ are assumed to be unknown, we can use a conservative minimization. It means that we minimize the nonrandom functional

(2.1) $$l(\theta, N) = \sum_{k=N+1}^{\infty} \theta_k^2 + V(N),$$

where $V(N)$ bounds from above the stochastic term $\sum_{k=1}^{N} \sigma_k^2 \xi_k^2$. It seems natural to choose this function such that

$$(2.2) \qquad \mathbf{E}\sup_N \left[\sum_{k=1}^{N} \sigma_k^2 \xi_k^2 - V(N)\right] \leq 0,$$

since then we can easily control the risk of any projection estimator with a data-driven bandwidth $\widetilde{N}$,

$$(2.3) \qquad \mathbf{E}_\theta \|\tilde{\theta}(\widetilde{N}) - \theta\|^2 \leq \mathbf{E}_\theta l(\theta, \widetilde{N}).$$

This motivation leads to the following definition: *a nonrandom function $\ell(\theta, N)$ such that $\mathbf{E}_\theta \sup_N [r(y, N) - \ell(\theta, N)] \leq 0$ is called a risk hull.*

Thus, we can say that $l(\theta, N)$ defined by (2.1) and (2.2) is a risk hull. Evidently, we want to have the upper bound (2.3) as small as possible. So, we are looking for the minimal hull. Note that this hull strongly depends on $\sigma_k^2$ and we present in the sequel a numerical recipe to compute it.

Once $V(N)$ satisfying (2.2) has been chosen, the minimization of $l(\theta, N)$ can be completed in the standard way by using unbiased estimation. Note that our problem is reduced to minimization of $-\sum_{k=1}^{N} \theta_k^2 + V(N)$. Replacing the unknown $\theta_k^2$ by their unbiased estimates $y_k^2 - \sigma_k^2$, we arrive at the following method of adaptive bandwidth choice:

$$\bar{N} = \arg\min_N \left[-\sum_{k=1}^{N} y_k^2 + \sum_{k=1}^{N} \sigma_k^2 + V(N)\right].$$

A cornerstone idea of this approach is that we can find a function $V(N)$ such that the data-driven $\bar{N}$ minimizes the risk hull $l(\theta, N)$ without significant losses, that is,

$$\mathbf{E}_\theta \bar{l}(\theta, \bar{N}) \lesssim \min_N l(\theta, N) + small\_term.$$

Therefore, combining this with (2.3), we get the inequality

$$(2.4) \qquad \mathbf{E}_\theta \|\tilde{\theta}(\bar{N}) - \theta\|^2 \lesssim \min_N l(\theta, N) + small\_term,$$

which represents a heuristic version of an oracle inequality for the RHM method.

Notice that when the risk is measured by the $l_2$-norm, RHM coincides with the empirical risk minimization approach which is usually used in model selection [4]. The major issue of model selection is the choice of a good penalization. In the framework of the RHM approach, this problem can be rephrased as follows: *to find the minimal risk hull, which can be minimized based on the data.* We do not believe that there is a good general formula for the optimal risk hull or for the penalty. What we can really do is to make use of the Monte Carlo method to compute an approximation of this hull. The goal of the present paper is to demonstrate that this approach works well for the regularization by projections.

2.2. *Statistical model and assumptions.* In the sequence space model (1.3), we supposed that $\sigma_k^2$ is a polynomially increasing sequence with $\sigma_1^2 > 0$. To be more precise, it is assumed that this sequence satisfies the following hypothesis.

POLYNOMIAL HYPOTHESIS. There exist constants $C_1, C_2, C_3$ such that for some $\beta \geq 0$ and for all $k > 1$

$$(2.5) \qquad C_1 \left( \frac{1}{2k} \sum_{i=1}^{2k} \sigma_i^4 \right)^{1/2} \leq \sigma_k^2 \leq C_2 \sigma_1^2 \left( \frac{1}{\sigma_1^2} \sum_{i=1}^{k-1} \sigma_i^2 \right)^{2\beta/(2\beta+1)}.$$

For any integer $s > 1$

$$(2.6) \qquad \frac{1}{\sigma_1^{2s}} \sum_{i=1}^{k} \sigma_i^{2s} \leq C_3^s \left( \frac{1}{\sigma_1^2} \sum_{i=1}^{k} \sigma_i^2 \right)^{(2s\beta+1)/(2\beta+1)}.$$

Let us comment very briefly on these assumptions. Assumption (2.5) means that $\sigma_k^2$ can have only polynomial growth. Indeed, since $x^{1/(2\beta+1)}$ is a concave function, we have by (2.5)

$$\left( \frac{1}{\sigma_1^2} \sum_{i=1}^{k} \sigma_i^2 \right)^{1/(2\beta+1)} - \left( \frac{1}{\sigma_1^2} \sum_{i=1}^{k-1} \sigma_i^2 \right)^{1/(2\beta+1)} \leq \frac{C_2}{2\beta+1},$$

and summing up these formulas, one can easily check that

$$\frac{1}{\sigma_1^2} \sum_{i=1}^{N} \sigma_i^2 \leq \left( \frac{C_2(N-1)}{2\beta+1} \right)^{2\beta+1} + 1,$$

(2.7)

$$\sigma_k^2 \leq C_2 \sigma_1^2 \left( \frac{C_2(k-1)}{2\beta+1} \right)^{2\beta} + C_2 \sigma_1^2.$$

Thus $\sigma_k$ can have only polynomial growth of order $\beta$, which we will call the *degree of inverse problem*.

2.3. *A risk hull.* The main ingredient of RHM is the function $U_0(k), k = 1,\ldots$, defined by (1.11). The simplest way to compute it is to make use of the Monte Carlo method. It should be mentioned that this method is time consuming since this function is related to large deviations of $\eta_k$. Lemma 1 below gives an asymptotic approximation for $U_0(k)$, but we will see that this approximation is not good for small $k$. Therefore we prefer to use the nonasymptotic formula (1.11) in our approach. It should be mentioned that the performance of RHM is sufficiently stable with respect to small perturbations of $U_0(k)$. Denote for brevity

$$\Sigma_N = \sum_{s=1}^{N} \sigma_s^4 \quad \text{and} \quad u_0(N) = \frac{U_0(N)}{\sqrt{2\Sigma_N}}.$$

LEMMA 1. *There exists an integer $N_0 \geq 1$ such that for all $N \geq N_0$*

$$(2.8) \qquad u_0(N) \geq u_1(N) \stackrel{\text{def}}{=} \sqrt{\log(\Sigma_N/(2\pi\sigma_1^4))}.$$

This fact plays a principal role in the proof of the following theorem.

THEOREM 2. *There exists a constant $C_*$ such that for any $\alpha > 0$*

$$(2.9) \qquad l_{\text{rhm}}(\theta, N) = \sum_{k=N+1}^{\infty} \theta_k^2 + \sum_{k=1}^{N} \sigma_k^2 + (1+\alpha)U_0(N) + \frac{C_* \sigma_1^2}{\alpha}$$

*is a risk hull, that is, $\mathbf{E} \sup_N [r(y, N) - l_{\text{rhm}}(\theta, N)] \leq 0$.*

This theorem says that uniformly in $N$, the loss $r(y, N)$ can be bounded by the risk hull $l_{\text{rhm}}(\theta, N)$. Thus, for any $\tilde{N}$ data-dependent, we can bound the risk of the projection regularization method by the expectation of the risk hull [see (2.3)].

We have mentioned that the URE and RHM methods can be viewed as minimizers of the penalized empirical risk [see (1.9)]. While the penalty corresponding to the URE is given by $\text{pen}_{\text{ure}}(N) = \sum_{k=1}^{N} \sigma_k^2$, the RHM method has the larger penalty $\text{pen}_{\text{rhm}}(N) = \sum_{k=1}^{N} \sigma_k^2 + (1+\alpha)U_0(N)$. Thus, it would be instructive to look at the ratio $\text{pen}_{\text{rhm}}(N)/\text{pen}_{\text{ure}}(N)$. If we suppose for a moment that the distribution $\eta_k$ can be approximated by a Gaussian law, then we get from (1.11)

$$U_0(N) \approx \widetilde{U}_0(N) = \sqrt{2\Sigma_N \log[\Sigma_N/(\pi\sigma_1^4)]}.$$

Under the polynomial hypothesis [see (2.5), (2.6)] it is easy to check that

$$U_0(N) = o\left(\sum_{k=1}^{N} \sigma_k^2\right), \qquad \widetilde{U}_0(N) = o\left(\sum_{k=1}^{N} \sigma_k^2\right), \qquad N \to \infty.$$

Nevertheless it is instructive to look at what is going on when $N$ is small. Therefore we plotted in Figure 2 the functions

$$\rho(N) = \frac{\text{pen}_{\text{rhm}}(N)}{\text{pen}_{\text{ure}}(N)} = 1 + \frac{(1+\alpha)U_0(N)}{\sum_{k=1}^{N} \sigma_k^2}, \qquad \widetilde{\rho}(N) = 1 + \frac{(1+\alpha)\widetilde{U}_0(N)}{\sum_{k=1}^{N} \sigma_k^2},$$

with $\alpha = 0.1$. Since we used the Monte Carlo method, the function $\rho(N)$ looks a little bit wiggly. The first display (direct estimation) shows that $(1+\alpha)U_0(N)$ is smaller than $\sum_{k=1}^{N} \sigma_k^2$ and this function cannot substantially affect the performance of URE. On the other hand, the second plot distinctly demonstrates that $(1+\alpha)U_0(N)$ dominates $\sum_{k=1}^{N} \sigma_k^2$ when $\sigma_k = \varepsilon k$. It means that in this case RHM and URE may work quite differently. Note also that in the case of inverse estimation the difference between $U_0(N)$ and its Gaussian approximation $\widetilde{U}_0(N)$ may be significant for small $N$. Certainly, $\widetilde{U}_0(N)/U_0(N) \to 1$ as $N \to \infty$, but very often numerical performance of RHM strongly depends on the behavior of the penalty function for small $N$, and this is why we used $U_0(N)$ in our method.

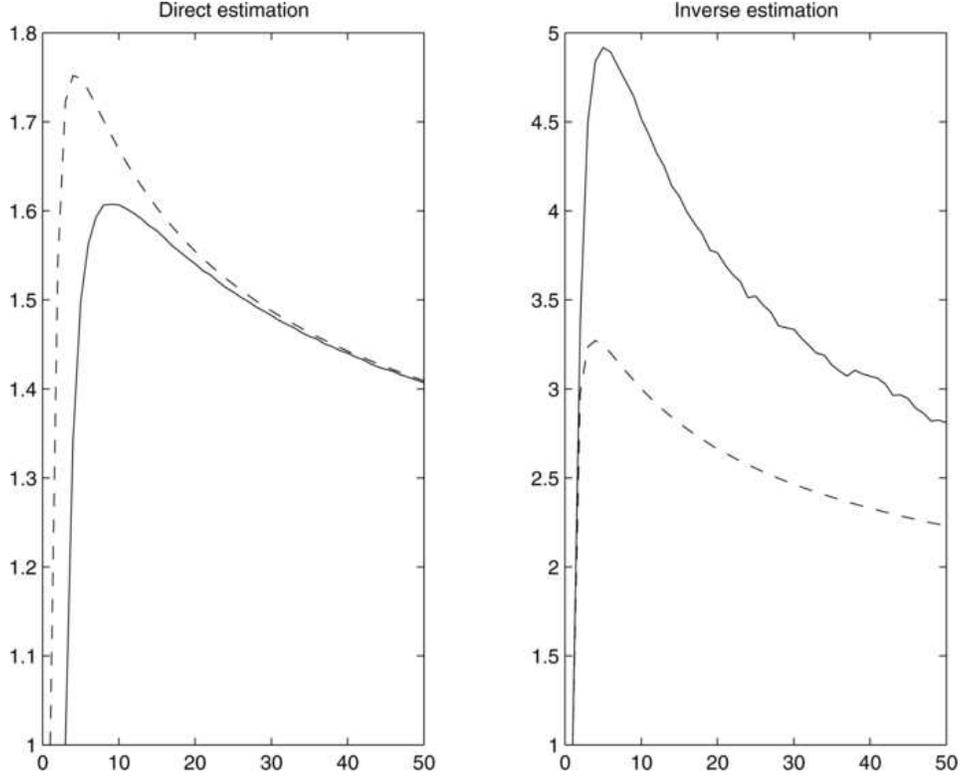

FIG. 2. *The functions $\rho(N)$ (solid line) and $\widetilde{\rho}(N)$ (dashed line) for direct ($\sigma_k = \varepsilon$) and inverse ($\sigma_k = \varepsilon k$) estimation.*

2.4. *The risk hull approach and URE.* Let us finish this section with a discussion of the URE method, which can be also viewed as a risk hull method. The following theorem justifies this idea.

THEOREM 3. *There exists a constant $C_u$ such that for any $\alpha > 0$*

$$l_{\mathrm{ure}}(\theta, N) = (1+\alpha)\left[\sum_{k=N+1}^{\infty} \theta_k^2 + \sum_{k=1}^{N} \sigma_k^2\right] + \frac{C_u}{\alpha^{4\beta+1}}\sigma_1^2$$

*is a risk hull.*

It is clear that the data-driven bandwidth choice $N_{\text{ure}}$ defined by (1.7) can be viewed as the minimization of the risk hull $l_{\text{ure}}(\theta, N)$. The following theorem provides an upper bound for the risk of this method.

THEOREM 4. *There exist constants $C^* > 0$ and $\gamma_0 > 0$ such that for all $\gamma \in (0, \gamma_0]$*

$$(2.10) \qquad \mathbf{E}_\theta \|\tilde{\theta}(N_{\text{ure}}) - \theta\|^2 \leq (1+\gamma) \inf_N R(\theta, N) + \frac{C^* \sigma_1^2}{\gamma^{4\beta+1}}.$$

This result rectifies Theorem 1 in [5]. It shows in particular that there is no logarithmic factor in the corresponding oracle inequality. At the first glance it seems that URE method may work better than RHM. This naive idea is motivated by the fact that

$$\inf_N R(\theta, N) < \inf_N R_{\text{rhm}}(\theta, N).$$

Recall that the left-hand side of this display represents the main term of the upper bound (2.10) while the right-hand side is the principal term of (1.12). But the real situation is not so trivial. In order to compare the bounds (2.10) and (1.12), we should take into account the remainder terms defined by constants $C_*$ and $C^*$. Both these constants depend on $\beta$ but their statistical nature and behavior are quite different, which follows from inspection of the proofs of Theorems 1 and 4. The constant $C^*$ may be very large even for $\beta > 1$ whereas $C_*$ remains moderate. We shall clearly see this phenomenon in the following section devoted to numerical simulations, but now let us discuss at the heuristic level the principal difficulties of URE. The basic idea of this method is that $\bar{R}(y, N) = -\sum_{k=1}^N y_k^2 + 2\sum_{k=1}^N \sigma_k^2$ is a good estimator for $\mathbf{E}_\theta \bar{R}(y, N) = -\sum_{k=1}^N \theta_k^2 + \sum_{k=1}^N \sigma_k^2$. In order to see that this idea may fail, it suffices to look at the variance

$$\mathbf{E}_\theta [\bar{R}(y, N) - \mathbf{E}_\theta \bar{R}(y, N)]^2 \geq 2 \sum_{k=1}^N \sigma_k^4.$$

So, $\bar{R}(y, N)$ might be considered a good estimator, if

$$\mathbf{E}_\theta \bar{R}(y, N) > 2 \left( 2 \sum_{k=1}^N \sigma_k^4 \right)^{1/2}.$$

This entails, in particular, that the following inequality should hold:

$$(2.11) \qquad \sum_{k=1}^N \sigma_k^2 > 2 \left( 2 \sum_{k=1}^N \sigma_k^4 \right)^{1/2},$$

for all $N \geq 1$. Notice that the factor 2 in the above inequality is, in some sense, very optimistic. In fact, it should be replaced by a function which tends to infinity as $N \to \infty$. However, let us suppose that $\sigma_k = \varepsilon k^\beta$ and look for integers $N_\beta$ for which (2.11) starts to work. For $\beta = 0$, we get $N_0 = 8$, for $\beta = 1$, $N_1 = 14$ and so on. It is easy to see that URE will always choose a bandwidth of order at least $N_\beta$. This evidently results in the risk order $\varepsilon^2 N_\beta^{2\beta+1}$. We would like to draw attention to the fact that this lower bound does not depend on the risk of the oracle $\inf_N R(\theta, N)$. The latter may be small while $\varepsilon^2 N_\beta^{2\beta+1}$ is large. Thus, roughly speaking, URE works well when

$$\inf_N R(\theta, N) > \varepsilon^2 N_\beta^{2\beta+1}.$$

Otherwise it fails. Unfortunately, the factor $N_\beta^{2\beta+1}$ is large even for moderate $\beta$; for $\beta = 1$ it is of order $10^3$.

The second almost evident fact is that the bandwidth $N$ of the best projection method is typically small when we deal with ill-posed problems. For instance, consider the minimax recovering of vectors $\theta$ from the Sobolev ball

$$W_m(L) = \left\{\theta : \sum_{k=1}^\infty \theta_k^2 k^{2m} \leq L\right\}.$$

Then it is easy to see (see, for details, [5]) that $N$ is of order $\varepsilon^{-2/(2m+2\beta+1)}$. Thus, when $\beta = 1$ and $m = 1$, this term is of order $\varepsilon^{-2/5}$. Therefore even for a very small noise level $\varepsilon^2 = 10^{-6}$, $N$ will not be larger than 20. Combining this with the previous remark, we see that in this case URE may not work properly. From an asymptotic viewpoint everything goes smoothly, but unfortunately asymptotic arguments start to work for very small $\varepsilon$.

**3. Simulations.** In this section we present some numerical properties of the RHM approach. Numerical testing of nonparametric statistical methods is a very difficult and delicate problem. The goal of this section is very modest. We would like to illustrate graphically Theorems 1 and 4. To do that, we propose to measure statistical performance of a method $\tilde{N}$ by *oracle efficiency* defined by

$$e_{\mathrm{or}}(\theta, \tilde{N}) = \frac{\inf_N \mathbf{E}_\theta \|\tilde{\theta}(N) - \theta\|^2}{\mathbf{E}_\theta \|\tilde{\theta}(\tilde{N}) - \theta\|^2}.$$

It should be mentioned that we use the inverse of the ratio $r(\theta)$ from (1.8) since we want to get a good graphical representation of the performance. We have seen in the Introduction that $r(\theta)$ may vary from 1 to 2000 for the URE method. This results in a degenerate plot of $r(\theta)$. Therefore, in order to avoid this effect, we use $e_{\mathrm{or}}(\theta, \tilde{N})$ instead of $r(\theta)$.

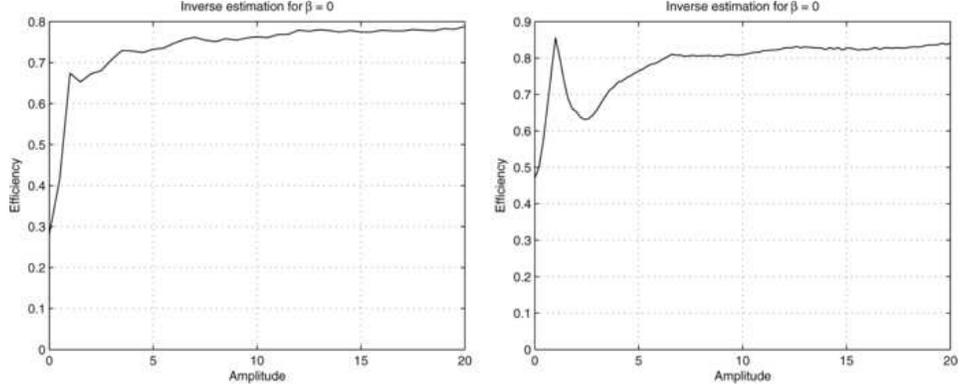

Fig. 3. *Oracle efficiency of URE and of RHM for direct estimation ($\sigma_k = \varepsilon$).*

Since it is evidently impossible to compute the oracle efficiency for all $\theta \in l_2$, we choose a sufficiently representative family of vectors $\theta$. In what follows we will use the linear family

$$\theta_i^a = \frac{a\varepsilon}{1 + (i/W)^m},$$

where $a$ defines amplitude, $W$ bandwidth and $m$ smoothness.

We shall vary $a$ in a large range and plot $r_{\mathrm{or}}(\theta^a, \tilde{N})$ as a function of $a$ which is directly related to the signal-to-noise ratio in the considered model. The parameters $m = 6$ and $W = 6$ are fixed. In other examples of $(W, m)$ the authors looked at, simulations showed that the oracle efficiency exhibits similar behavior.

Two methods of data-driven bandwidth choice will be compared: URE and RHM with $\alpha = 1.1$. It is easy to see that for these methods $r_{\mathrm{or}}(\theta^a, \tilde{N})$ does not depend on $\varepsilon$. This function was computed by the Monte Carlo method with 40,000 replications. We start with direct estimation where $\sigma_k \equiv \varepsilon$. Figure 3 shows the oracle efficiency of URE (left panel) and the oracle efficiency of RHM (right panel). Comparing these plots, one can say that both methods work reasonably well. However, if we deal with an inverse problem such as derivative estimation, we can see a significant difference between these methods. The corresponding oracle efficiencies are plotted on the left and right panels of Figure 4. For small values of $a$ the performance of URE is very poor, whereas RHM demonstrates very stable behavior. For very large $a = 500$ the oracle efficiency of URE is of order 0.16, while RHM always has efficiency greater than 0.4. Figure 5 deals with the case when the inverse problem becomes really ill-posed ($\sigma_k = \varepsilon k^2$). In this situation URE fails completely. Its maximal oracle efficiency is of order $3 * 10^{-4}$. Nevertheless, RHM has a good efficiency (greater than 0.3). In the context of Theorem 1 and Theorem 4 this example shows that the constants $C_*$ and $C^*$ are quite different: while $C_*$ is small, $C^*$ is really large. Unfortunately, it means that the terms which are asymptotically small in Theorem 4 may easily dominate the oracle risk.

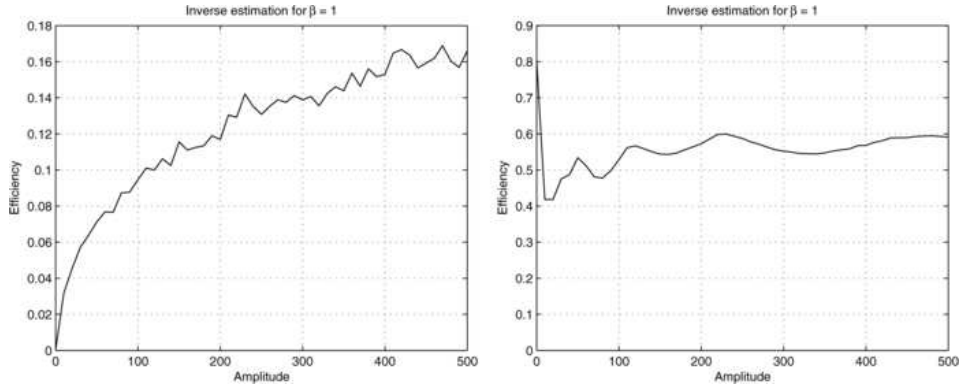

Fig. 4. *Oracle efficiency of URE and of RHM for first-order derivative estimation ($\sigma_k = \varepsilon k$).*

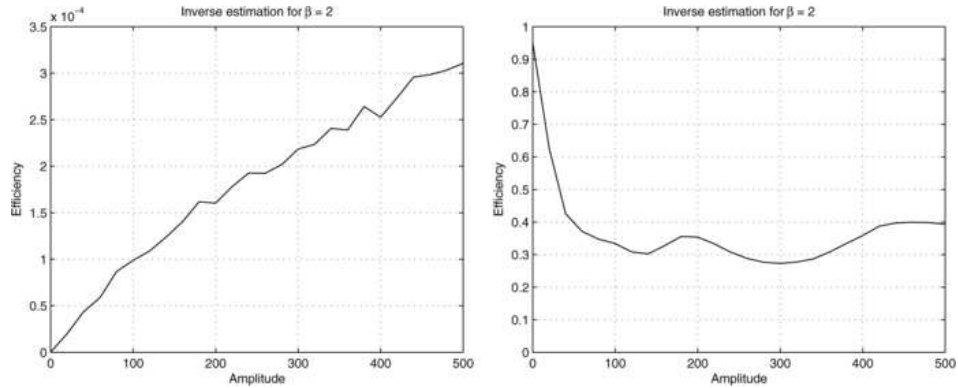

FIG. 5. *Oracle efficiency of URE and of RHM for second-order derivative estimation* $(\sigma_k = \varepsilon k^2)$.

Let us finish this section with a short discussion of the role played by $\alpha$. In the previous numerical simulations this parameter was 1.1. What happens if we set this parameter to 0? The answer depends on $\beta$. If $\beta$ is small, $\beta \leq 1$, everything goes smoothly. However, even for $\beta = 2$ this choice results in an instable procedure. On the other hand, taking $\alpha$ to be large leads to poor performance of RHM.

**4. Proofs.**

4.1. *Proof of Theorem* 2.

PROOF OF LEMMA 1. Denote for brevity
$$\kappa_N = \eta_N/\sqrt{2\Sigma_N}$$
and
$$\Phi_N(t) = \mathbf{E}\exp(it\kappa_N).$$

We begin with an upper bound for the absolute value of $\Phi_N(t)$. Recalling the definition of $\eta_N$ and using (2.5), we have

$$|\Phi_N(t)| \leq \exp\left[-\frac{1}{4}\sum_{l=1}^{N}\log\left(1+\frac{2t^2\sigma_l^4}{\Sigma_N}\right)\right]$$

$$\leq \exp\left[-\frac{N}{8}\log\left(1+\frac{2t^2\sigma_{N/2}^4}{\Sigma_N}\right)\right]$$

$$\leq \left(1+\frac{Ct^2}{N}\right)^{-N/8}.$$

With this inequality, we have that for all $x \leq \sqrt{N/C}$

$$\int_{|t|\geq x}|\Phi_N(t)|\,dt \leq \int_{x\leq|t|\leq\sqrt{2N/C}}|\Phi_N(t)|\,dt + \int_{|t|\geq\sqrt{2N/C}}|\Phi_N(t)|\,dt$$

(4.1)
$$\leq \int_{|t|\geq x}\exp[-Ct^2]\,dt + \int_{|t|\geq\sqrt{2N/C}}\left(\frac{Ct^2}{N}\right)^{-N/8}dt$$

$$\leq \exp[-Cx^2] + \sqrt{\frac{N^3}{8C}}2^{-N/8} \leq \exp[-Cx^2].$$

Let us fix an integer $M$. Then by the Taylor formula we get that for all $|t| \leq \sqrt{N/C}$

$$\Phi_N(t) = \exp\left\{-\frac{t^2}{2} + \sum_{s=3}^{M-1}\frac{(-i)^s 2^{s/2}R_s t^s}{s} + O\left(\frac{CR_M 2^{M/2}t^M}{M}\right)\right\},$$

where $R_s = (\Sigma_N)^{-s/2}\sum_{l=1}^{N}\sigma_l^{2s}$.

It follows easily from (2.5) that $\sigma_N^4 \leq C\sigma_{N/2}^4$. This gives $|R_s| \asymp N^{-s/2+1}$. Therefore, expanding $\Phi_N(t)\exp(t^2/2)$ into Taylor series, it is easy to see that there exist functions $Q_M(s,N), s=3,\ldots,M$, uniformly bounded in $N$ and $s$ such that

(4.2) $\Phi_N(t)\exp(t^2/2) = 1 + N\sum_{s=3}^{M-1}Q_M(s,N)\left(\frac{it}{\sqrt{N}}\right)^s + O\left(N\left(\frac{t}{\sqrt{N}}\right)^M\right).$

Define now the following approximation of $\Phi_N(t)$:

$$\Phi_N^M(t) = \exp(-t^2/2)\left[1 + N\sum_{s=3}^{M-1} Q_M(s,N)\left(\frac{it}{\sqrt{N}}\right)^s\right].$$

Now we can approximate the probability $\mathbf{P}(\kappa_N > x)$ by

$$(4.3) \quad P_N^M(x) = \int_x^\infty p_N^M(v)\,dv \qquad \text{with } p_N^M(v) = \frac{1}{2\pi}\int_{-\infty}^\infty \exp(-itv)\Phi_N^M(t)\,dt.$$

Notice that

$$(4.4) \quad P_N^M(x) = \phi(x) - \frac{k}{\sqrt{2\pi}}\sum_{s=3}^{M-1}(-1)^s Q_M(s,N)k^{-s/2}\frac{d^{s-1}}{dx^{s-1}}\exp(-x^2/2),$$

where

$$\phi(x) = \frac{1}{\sqrt{2\pi}}\int_x^\infty \exp(-u^2/2)\,du.$$

Then by the Parseval identity and (4.1), (4.2), we obtain

$$(4.5) \quad \begin{aligned} &|\mathbf{P}(\kappa_N > x) - P_N^M(x)| \\ &\leq \frac{1}{2\pi}\int_{-\infty}^\infty |t|^{-1}|\Phi_N^M(t) - \Phi_N(t)|\,dt \\ &= \int_{|t|\leq\sqrt{N/C}} \frac{|\Phi_N^M(t) - \Phi_N(t)|}{|t|}\,dt + \int_{|t|\geq\sqrt{N/C}} \frac{|\Phi_N^M(t) - \Phi_N(t)|}{|t|}\,dt \\ &\leq \frac{C}{N^{M/2}} + \exp[-CN] \leq \frac{C}{N^{M/2}}. \end{aligned}$$

Using (4.3) and (4.5), it is easy to see that

$$(4.6) \qquad \mathbf{P}(\kappa_N > x) \geq P_N^M(x) - \frac{C}{N^{M/2}}.$$

Now we are ready to complete the proof of the lemma. Since the function $F(x) = \mathbf{E}\kappa_N I(\kappa_N \geq x)$ is a monotone nondecreasing function in $x \geq 0$, we need to check that for sufficiently large $N$ [see (1.11) and (2.8)]

$$\mathbf{E}\kappa_N I(\kappa_N \geq u_1(N)) \geq \frac{\sigma_1^2}{\sqrt{2\Sigma_N}}.$$

It follows from the above equation and integration by parts that it suffices to show that

$$(4.7) \qquad u_1(N)\mathbf{P}(\kappa_N \geq u_1(N)) + \int_{u_1(N)}^{u_1(N)+1} \mathbf{P}(\kappa_N > x)\,dx \geq \frac{\sigma_1^2}{\sqrt{2\Sigma_N}}.$$

Using (4.6), we bound the left-hand side as

$$u_1(N)\mathbf{P}(\kappa_N \geq u_1(N)) + \int_{u_1(N)}^{u_1(N)+1} \mathbf{P}(\kappa_N > x)\,dx$$

$$\geq u_1(N)P_k^M(u_1(N)) + \int_{u_1(N)}^{u_1(N)+1} P_N^M(x)\,dx - \frac{Cu_1(N)}{N^{M/2}}$$

(4.8)

$$= u_1(N)\phi[u_1(N)] + \int_{u_1(N)}^\infty \phi(x)\,dx - \int_{u_1(N)+1}^\infty \phi(x)\,dx$$

$$- (1 + u_1(N)) \max_{u_1(N) \leq x \leq u_1(N)+1} |P_N^M(x) - \phi(x)| - \frac{Cu_1(N)}{N^{M/2}}.$$

Integrating by parts, we get

(4.9) $\quad u_1(N)\phi[u_1(N)] + \int_{u_1(N)}^\infty \phi(x)\,dx = \frac{1}{\sqrt{2\pi}} e^{-u_1^2(N)/2} = \frac{2\sigma_1^2}{\sqrt{2\Sigma_N}}.$

Noticing that in view of (2.7)

(4.10) $\quad\quad\quad \sqrt{\log(N/(2\pi))} \leq u_1(N) \leq C\sqrt{\log(N)}$

and integrating by parts, we have as $N \to \infty$

(4.11) $\displaystyle\int_{u_1(N)+1}^\infty \phi(x)\,dx \leq Ce^{-(u_1(N)+1)^2/2} \leq \frac{2\sigma_1^2}{\sqrt{2\Sigma_N}} e^{-u_1(N)} = o\left(\frac{\sigma_1^2}{\sqrt{2\Sigma_N}}\right),$

and by (4.4) and (4.10),

$$(1 + u_1(N)) \max_{u_1(N) \leq x \leq u_1(N)+1} |P_N^M(x) - \phi(x)|$$

(4.12)

$$\leq \frac{C\sigma_1^2 u_1^3(N)}{\sqrt{2N\Sigma_N}} = o\left(\frac{\sigma_1^2}{\sqrt{2\Sigma_N}}\right).$$

Finally, note that we can choose sufficiently large $M$ such that [see (4.10) and (2.7)]

$$\frac{Cu_1(N)}{N^{M/2}} = o\left(\frac{\sigma_1^2}{\sqrt{2\Sigma_N}}\right), \quad k \to \infty.$$

Combining this equation with (4.8)–(4.12) we arrive at (4.7), thus finishing the proof of the lemma. □

LEMMA 2. *For some $C > 0$*

(4.13) $\quad\quad \mathbf{P}\{\eta_N > x\} \leq \exp\left(-\frac{Cx^2}{\Sigma_N}\right), \quad 0 \leq x \leq \frac{\Sigma_N}{\sigma_N^2}.$

PROOF. Certainly, this fact is well known and we prove it only for the reader's convenience. We use the inequality $\log(y) \geq y - 1 - (1 - 1/y)^2/2, y \in (0, 1]$, which can be checked easily since the first derivative in $y$ of $\log(y) - y + 1 + (1 - 1/y)^2/2$ is negative. Therefore, for any positive $\lambda$

$$\mathbf{E}\exp(\lambda \eta_N) = \exp\left\{-\lambda \sum_{i=1}^N \sigma_i^2 - \frac{1}{2}\sum_{i=1}^N \log(1 - 2\lambda\sigma_i^2)_+\right\}$$

$$(4.14) \qquad \leq \exp\left\{\lambda^2 \sum_{i=1}^N \frac{\sigma_i^4}{(1 - 2\lambda\sigma_i^2)_+^2}\right\}$$

$$\leq \exp\left\{\lambda^2 \sum_{i=1}^N \sigma_i^4 + 4\lambda^3 \sum_{i=1}^N \frac{\sigma_i^6}{(1 - 2\lambda\sigma_i^2)_+^2}\right\}.$$

Then, by the Markov inequality, we have

$$\mathbf{P}\{\eta_N > x\} \leq \exp(-\lambda x)\mathbf{E}\exp(\lambda \eta_N)$$

$$\leq \exp\left\{-\lambda x + \lambda^2 \sum_{i=1}^N \frac{\sigma_i^4}{(1 - 2\lambda\sigma_i^2)_+^2}\right\}.$$

In order to prove (4.13), we take $\lambda = x/(8\Sigma_N)$. □

Define the auxiliary function

$$U(\alpha) = -1 - \frac{1}{2\alpha}\log(1 - 2\alpha), \qquad \alpha \in (0, 1/2).$$

By the Taylor formula $U(\alpha) = 2\alpha \sum_{i=2}^\infty (2\alpha)^{i-2}/i$. This yields immediately that $\alpha \leq U(\alpha) \leq \alpha/(1 - 2\alpha)$ and

$$(4.15) \qquad \frac{\alpha}{1 + 2\alpha} \leq U^{-1}(\alpha) \leq \alpha, \qquad \alpha > 0,$$

where $U^{-1}(\alpha)$ denotes the inverse function.

LEMMA 3. *Let $S_N = \sum_{i=1}^N b_i^2(\xi_i^2 - 1) - U(\alpha)\sum_{i=1}^N b_i^4$ where $\xi_i$ are i.i.d. $\mathcal{N}(0,1)$ and $b_i^2 \leq 1$. Then for any $\alpha \in (0, 1/2)$*

$$(4.16) \qquad \mathbf{E}\sup_{N \geq 1} S_N \leq \alpha^{-1}$$

*and*

$$(4.17) \qquad \mathbf{P}\left(\sup_{N \geq 1} S_N > x\right) \leq \exp(-\alpha x).$$

The proof follows from the Doob inequality (see, e.g., [13]).

PROOF OF THEOREM 2. Define the exponential grid

(4.18) $$n_s = \lfloor (1 + p\sqrt{\alpha})^s \rfloor,$$

where $p$ is a sufficiently small constant which will be chosen later on. By Lemma 1 and a simple algebra we have

$$\mathbf{E} \sup_N \{\eta_N - (1+\alpha)U_0(N)\}$$

(4.19)
$$\leq \sum_{s=1}^{\infty} \mathbf{E} \max_{n_s \leq N < n_{s+1}} [\eta_N - (1+\alpha)U_0(N)]_+$$

$$\leq \sum_{s=1}^{\infty} \mathbf{E} \max_{n_s \leq N < n_{s+1}} [\eta_N - (1+\alpha)\sqrt{2\Sigma_k \log(C\Sigma_{n_s}/\sigma_1^4)}]_+$$

$$\leq \sum_{s=1}^{\infty} \mathbf{E}[\eta_{n_s} - (1+\alpha)\sqrt{2\Sigma_{n_s} \log(C\Sigma_{n_s}/\sigma_1^4)} + \epsilon_s]_+,$$

where

$$\epsilon_s = \max_{n_s < N < n_{s+1}} \left\{ \sum_{i=n_s+1}^{N} \sigma_i^2(\xi_i^2 - 1) - (1+\alpha)\frac{\sqrt{2\log(C\Sigma_{n_s}/\sigma_1^4)}}{2\sqrt{\Sigma_{n_{s+1}}}}[\Sigma_N - \Sigma_{n_s}] \right\}.$$

Denote for brevity

$$\Gamma_s = (1+\alpha)\frac{\sqrt{2\log(C\Sigma_{n_s}/\sigma_1^4)}}{2\sqrt{\Sigma_{n_{s+1}}}} \quad \text{and} \quad A_s = (1+\alpha)\sqrt{\log(C\Sigma_{n_s}/\sigma_1^4)}.$$

Then by (4.15) and (4.17) we obtain

$$\mathbf{P}\{\epsilon_s \geq x\}$$

$$= \mathbf{P}\left\{ \max_{n_s < N < n_{s+1}} \left[ \sum_{i=n_s+1}^{N} \sigma_i^2(\xi_i^2 - 1) - \Gamma_s \sum_{i=n_s+1}^{N} \sigma_i^4 \right] \geq x \right\}$$

$$= \mathbf{P}\left\{ \max_{n_s < N < n_{s+1}} \left[ \sum_{i=n_s+1}^{N} \frac{\sigma_i^2}{\sigma_{n_{s+1}}^2}(\xi_i^2 - 1) - \Gamma_s \sigma_{n_{s+1}}^2 \sum_{i=n_s+1}^{N} \frac{\sigma_i^4}{\sigma_{n_{s+1}}^4} \right] \geq \frac{x}{\sigma_{n_{s+1}}^2} \right\}$$

$$\leq \exp[-U^{-1}(\Gamma_s \sigma_{n_{s+1}}^2)x/\sigma_{n_{s+1}}^2] \leq \exp[-\Gamma_s(1 - 2\Gamma_s \sigma_{n_{s+1}}^2)x],$$

or equivalently

(4.20) $$\mathbf{P}\left\{ \frac{\epsilon_s}{\sqrt{2\Sigma_{n_s}}} \geq x \right\} \leq \exp[-A_s(\Sigma_{n_s}/\Sigma_{n_{s+1}})^{1/2}(1 - 2\Gamma_s \sigma_{n_{s+1}}^2)x].$$

Notice that $1 - 2\Gamma_s \sigma_{n_{s+1}}^2 \geq 0$ for sufficiently large $s$. Using (4.20) and integrating by parts, we get

$$\mathbf{E}[\eta_{n_s} - (1+\alpha)\sqrt{2\Sigma_{n_s} \log(C\Sigma_{n_s}/\sigma_1^4)} + \epsilon_s]_+$$

$$= \sqrt{2\Sigma_{n_s}} \mathbf{E}\left[\frac{\eta_{n_s}}{\sqrt{2\Sigma_{n_s}}} + \frac{\epsilon_s}{\sqrt{2\Sigma_{n_s}}} - A_s\right]_+$$

$$= \sqrt{2\Sigma_{n_s}} \int_{A_s}^\infty \mathbf{P}\left\{\frac{\eta_{n_s}}{\sqrt{2\Sigma_{n_s}}} + \frac{\epsilon_s}{\sqrt{2\Sigma_{n_s}}} \geq x\right\} dx$$

(4.21)

$$\leq \sqrt{2\Sigma_{n_s}} \int_{A_s}^\infty \mathbf{E} \exp\left\{-A_s\sqrt{\frac{\Sigma_{n_s}}{\Sigma_{n_{s+1}}}}(1 - 2\Gamma_s\sigma_{n_{s+1}}^2)\left(x - \frac{\eta_{n_s}}{\sqrt{2\Sigma_{n_s}}}\right)\right\} dx$$

$$= \frac{\sqrt{2\Sigma_{n_{s+1}}}}{A_s} \exp\{-A_s^2(\Sigma_{n_s}/\Sigma_{n_{s+1}})^{1/2}(1 - 2\Gamma_s\sigma_{n_{s+1}}^2)\}$$

$$\times \mathbf{E} \exp\left\{A_s(1 - 2\Gamma_s\sigma_{n_{s+1}}^2)\frac{\eta_{n_s}}{\sqrt{2\Sigma_{n_{s+1}}}}\right\}.$$

In order to bound from above the last term in this inequality, we have by (4.14) that for any positive $\lambda$

$$\mathbf{E} \exp\{\lambda \eta_{n_s}\} \leq \exp\left\{\lambda^2 \sum_{i=1}^{n_s} \sigma_i^4 + 4\lambda^3 \sum_{i=1}^{n_s} \frac{\sigma_i^6}{(1 - 2\lambda\sigma_i^2)_+^2}\right\}.$$

Using this inequality with $\lambda = A_s(1 - 2\Gamma_s\sigma_{n_{s+1}}^2)/\sqrt{2\Sigma_{n_{s+1}}}$ and noticing that by virtue of the polynomial hypothesis $\lambda^3 \sum_{i=1}^{n_s} \sigma_i^6 (1 - 2\lambda\sigma_i^2)_+^{-2} \leq C$, we immediately get

$$\mathbf{E} \exp\left\{A_s \frac{\eta_{n_s}}{\sqrt{2\Sigma_{n_{s+1}}}}(1 - 2\Gamma_s\sigma_{n_{s+1}}^2)\right\} \leq C \exp\left\{\frac{A_s^2 \Sigma_{n_s}}{2\Sigma_{n_{s+1}}}(1 - 2\Gamma_s\sigma_{n_{s+1}}^2)\right\}.$$

Therefore, combining this with (4.21), we obtain

$$\mathbf{E}[\eta_{n_s} - (1+\alpha)\sqrt{2\Sigma_{n_s} \log(C\Sigma_{n_s}/\sigma_1^4)} + \epsilon_s]_+$$

(4.22)

$$\leq C \frac{\sqrt{2\Sigma_{n_{s+1}}}}{A_s} \exp\left\{-\frac{A_s^2}{2}(1 - 2\Gamma_s\sigma_{n_{s+1}}^2)\left(2\sqrt{\frac{\Sigma_{n_s}}{\Sigma_{n_{s+1}}}} - \frac{\Sigma_{n_s}}{\Sigma_{n_{s+1}}}\right)\right\}$$

$$\leq C \frac{\sqrt{2\Sigma_{n_{s+1}}}}{A_s} \exp\left\{-\frac{A_s^2}{2}(1 - 2\Gamma_s\sigma_{n_{s+1}}^2)\left[1 - \frac{1}{4}\left(\frac{\Sigma_{n_s} - \Sigma_{n_{s+1}}}{\Sigma_{n_{s+1}}}\right)^2\right]\right\}.$$

Let us choose now the parameter $p$ of the exponential grid. Note that by (2.5)

$$\left(\frac{\Sigma_{n_s} - \Sigma_{n_{s+1}}}{\Sigma_{n_{s+1}}}\right)^2 \leq \frac{C\sigma_{n_{s+1}}^8(n_{s+1} - n_s)^2}{\sigma_{n_{s+1}}^8 n_{s+1}^2} \leq \frac{Cp^2\alpha}{(1 - p\sqrt{\overline{\alpha}})^2}.$$

Thus it is clear that we can always choose a sufficiently small $p$ such that

$$\left(\frac{\Sigma_{n_s} - \Sigma_{n_{s+1}}}{\Sigma_{n_{s+1}}}\right)^2 \leq 4\alpha.$$

Hence from (4.22) we get

$$\mathbf{E}[\eta_{n_s} - (1+\alpha)\sqrt{2\Sigma_{n_s}\log(C\Sigma_{n_s}/\sigma_1^4)} + \epsilon_s]_+$$

(4.23)
$$\leq CA_s^{-1}\sqrt{\Sigma_{n_{s+1}}}\exp[-(1+\alpha-2\alpha^2)(1-2\Gamma_s\sigma_{n_{s+1}}^2)\log(\sqrt{C\Sigma_{n_s}/\sigma_1^4})]$$
$$\leq C\sigma_1^2 A_s^{-1}\exp[-\alpha(1-2\alpha)\log(\sqrt{C\Sigma_{n_s}/\sigma_1^4})]$$
$$\leq C\sigma_1^2 s^{-1/2}\alpha^{-1/4}\exp(-Cp\alpha^{3/2}s).$$

In the above inequality we used the fact that $A_s \geq (1+\alpha)\sqrt{\log(n_s)}$ and that $\Gamma_s \sigma_{n_{s+1}}^2 \log(C\Sigma_{n_s}/\sigma_1^4)$ is uniformly bounded in $s$. Finally, substituting (4.23) in (4.19), we have

$$(4.24) \quad \mathbf{E}\sup_N\{\eta_N - (1+\alpha)U_0(N)\} \leq \frac{C\sigma_1^2}{\sqrt[4]{\alpha}}\sum_{s=1}^\infty \frac{\exp(-Cp\alpha^{3/2}s)}{\sqrt{s}} \leq \frac{C\sigma_1^2}{\alpha\sqrt{p}},$$

thus proving the theorem. $\square$

4.2. *Proofs of Theorems* 1, 3 *and* 4. We start with two technical lemmas. Their proofs can be found in [13].

LEMMA 4. *Let $\kappa \geq 1$ be an integer random variable. Then for any $N = 1, 2, \ldots$*

$$\mathbf{E}\sum_{i=\kappa}^\infty \sigma_i\theta_i\xi_i \geq -\left\{3\sigma_N^2\mathbf{E}\sum_{i=\kappa}^\infty \theta_i^2 + 3\mathbf{E}\sigma_\kappa^2\sum_{i=N}^\infty \theta_i^2\right\}^{1/2}.$$

LEMMA 5. *For any $Q \in (1/2, (2\beta+1)/(4\beta+1)]$ there exist constants $C(Q) > 0$ and $\alpha(Q) > 0$ such that for all $\alpha \in (0, \alpha(Q))$ the following inequality holds:*

$$(4.25) \quad \mathbf{E}\sup_{N\geq 1}\left\{\frac{\eta_N}{\sigma_1^2} - \alpha\left(\frac{\Sigma_N}{\sigma_1^4}\right)^Q\right\} \leq C(Q)\alpha^{-1/(2Q-1)}.$$

PROOF OF THEOREM 1. In view of Theorem 2, for any $\mu > 0$

$$l_\mu(\theta, N) = \sum_{i=N+1}^\infty \theta_i^2 + \sum_{i=1}^N \sigma_i^2 + (1+\mu)U_0(N) + \frac{C\sigma_1^2}{\mu}$$

is a risk hull, and therefore

(4.26) $$\mathbf{E}_\theta \|\tilde{\theta}(N_{\mathrm{rhm}}) - \theta\|^2 \leq \mathbf{E}_\theta l_\mu(\theta, N_{\mathrm{rhm}}).$$

On the other hand, since $N_{\mathrm{rhm}}$ minimizes $\bar{R}_{\mathrm{pen}}(y, N)$ [see (1.9)], we have for any integer $N$

(4.27) $$\mathbf{E}_\theta \bar{R}_{\mathrm{pen}}(y, N_{\mathrm{rhm}}) \leq \mathbf{E}_\theta \bar{R}_{\mathrm{pen}}(y, N) = R_{\mathrm{rhm}}(\theta, N) + \|\theta\|^2.$$

In order to combine the inequalities (4.26) and (4.27), we rewrite $l_\mu(\theta, N_{\mathrm{rhm}})$ in terms of $\bar{R}_{\mathrm{pen}}(y, N_{\mathrm{rhm}})$,

$$\bar{R}_{\mathrm{pen}}(y, N_{\mathrm{rhm}}) + \|\theta\|^2 + \frac{C\sigma_1^2}{\mu}$$
$$= l_\mu(\theta, N_{\mathrm{rhm}}) - 2 \sum_{i=1}^{N_{\mathrm{rhm}}} \sigma_i \theta_i \xi_i - \sum_{i=1}^{N_{\mathrm{rhm}}} \sigma_i^2(\xi_i^2 - 1) + (\alpha - \mu) U_0(N_{\mathrm{rhm}}).$$

Therefore, using this equation and (4.26), (4.27), we obtain that for any integer $N$

(4.28) $$\mathbf{E}_\theta \|\tilde{\theta}(N_{\mathrm{rhm}}) - \theta\|^2 \leq R_{\mathrm{rhm}}(\theta, N) + \frac{C\sigma_1^2}{\mu} + 2\mathbf{E}_\theta \sum_{i=1}^{N_{\mathrm{rhm}}} \sigma_i \theta_i \xi_i$$
$$+ \mathbf{E}_\theta \left[ \sum_{i=1}^{N_{\mathrm{rhm}}} \sigma_i^2(\xi_i^2 - 1) - (\alpha - \mu) U_0(N_{\mathrm{rhm}}) \right].$$

Our next step is to control the last two terms in the above equation. By Lemma 4 we have that for any driven bandwidth $\tilde{N}$

(4.29) $$\mathbf{E}_\theta \sum_{i=1}^{\tilde{N}} \sigma_i \theta_i \xi_i = -\mathbf{E}_\theta \sum_{i=\tilde{N}+1}^{\infty} \sigma_i \theta_i \xi_i$$
$$\leq 2|\sigma_N| \left( \mathbf{E}_\theta \sum_{i=\tilde{N}+1}^{\infty} \theta_i^2 \right)^{1/2} + 2 \left( \sum_{i=N+1}^{\infty} \theta_i^2 \right)^{1/2} \sqrt{\mathbf{E}_\theta \sigma_{\tilde{N}}^2}.$$

Noticing that by (2.5) $\sigma_k^2 \leq C\sigma_1^2(\sigma_1^{-2} \sum_{i=1}^k \sigma_i^2)^{2\beta/(2\beta+1)}$ and using the Young inequality,

(4.30) $$xy^r \leq ry + (1-r)x^{1/(1-r)}, \qquad r \in (0,1),$$

with $r = 1/2$ and (4.29), we get that for any $\gamma > 0$

$$\mathbf{E}_\theta \sum_{i=1}^{\tilde{N}} \sigma_i \theta_i \xi_i$$

$$\leq C|\sigma_1| \left(\frac{1}{\sigma_1^2}\sum_{i=1}^{N}\sigma_k^2\right)^{\beta/(2\beta+1)} \left(\mathbf{E}_\theta \sum_{i=\tilde{N}+1}^{\infty}\theta_i^2\right)^{1/2}$$

$$+ C|\sigma_1|\left(\mathbf{E}_\theta \frac{1}{\sigma_1^2}\sum_{i=1}^{\tilde{N}}\sigma_k^2\right)^{\beta/(2\beta+1)} \left(\sum_{i=N+1}^{\infty}\theta_i^2\right)^{1/2}$$

$$\leq \gamma \sum_{i=N+1}^{\infty}\theta_i^2 + \gamma \mathbf{E}_\theta \sum_{i=\tilde{N}+1}^{\infty}\theta_i^2$$

$$+ \frac{C\sigma_1^{2/(2\beta+1)}}{\gamma^{(4\beta+1)/(2\beta+1)}}\left[\left(\gamma\sum_{i=1}^{N}\sigma_k^2\right)^{2\beta/(2\beta+1)} + \left(\gamma\mathbf{E}_\theta\sum_{i=1}^{\tilde{N}}\sigma_k^2\right)^{2\beta/(2\beta+1)}\right].$$

Once again using (4.30) with $r = 2\beta/(2\beta+1)$, we continue the above inequality as follows:

$$\mathbf{E}_\theta \sum_{i=1}^{\tilde{N}}\sigma_i\theta_i\xi_i$$

(4.31)
$$\leq \gamma\left(\sum_{i=N+1}^{\infty}\theta_i^2 + \sum_{i=1}^{N}\sigma_k^2\right) + \gamma\mathbf{E}_\theta\left(\sum_{i=\tilde{N}+1}^{\infty}\theta_i^2 + \sum_{i=1}^{\tilde{N}}\sigma_k^2\right) + \frac{C\sigma_1^2}{\gamma^{4\beta+1}}$$

$$\leq \gamma R(\theta, N) + \gamma \mathbf{E}_\theta \|\tilde{\theta}(\tilde{N}) - \theta\|^2 - \gamma \mathbf{E}_\theta \sum_{i=1}^{\tilde{N}}\sigma_k^2(\xi_i^2 - 1) + \frac{C\sigma_1^2}{\gamma^{4\beta+1}}.$$

Therefore, substituting (4.31) in (4.28) and then using (4.24), we obtain

$$(1-\gamma)\mathbf{E}_\theta\|\tilde{\theta}(N_{\text{rhm}}) - \theta\|^2$$

$$\leq (1+\gamma)R_{\text{rhm}}(\theta, N) + \frac{C\sigma_1^2}{\mu} + \frac{C\sigma_1^2}{\gamma^{4\beta+1}}$$

$$+ (1-\gamma)\mathbf{E}_\theta\left[\sum_{i=1}^{N_{\text{rhm}}}\sigma_i^2(\xi_i^2 - 1) - \frac{\alpha-\mu}{1-\gamma}U_0(N_{\text{rhm}})\right]$$

$$\leq (1+\gamma)R_{\text{rhm}}(\theta, N) + \frac{C\sigma_1^2}{\mu} + \frac{C\sigma_1^2}{\gamma^{4\beta+1}} + \frac{(1-\gamma)^2 C\sigma_1^2}{(\alpha-\mu+\gamma-1)_+}.$$

Finally, choosing $\mu = \gamma$, completes the proof. $\square$

PROOF OF THEOREM 3. This suffices to show that for any sufficiently small $\alpha > 0$

$$\mathbf{E}\sup_k\left[\eta_k - \alpha\sum_{i=1}^{k}\sigma_i^4\right] \leq \frac{C_u}{\alpha^{4\beta+1}}\sigma_1^2.$$

In view of (2.6) the proof follows immediately from Lemma 5 with $Q = (2\beta+1)/(4\beta+1)$. $\square$

PROOF OF THEOREM 4. This follows the main lines of the proof of Theorem 1. By Theorem 3 we have

$$(4.32) \quad \mathbf{E}_\theta \|\tilde{\theta}(N_{\text{ure}}) - \theta\|^2 \leq \mathbf{E}_\theta l_{\text{ure}}(\theta, N_{\text{ure}}) = (1+\alpha)\mathbf{E}_\theta R(\theta, N_{\text{ure}}) + \frac{C_u}{\alpha^{4\beta+1}}\sigma_1^2.$$

Since $N_{\text{ure}}$ minimizes $-\sum_{i=1}^N y_i^2 + 2\sum_{i=1}^N \sigma_i^2$, we get for any integer $N$

$$(4.33) \qquad -\sum_{i=1}^{N_{\text{ure}}} y_i^2 + 2\sum_{i=1}^{N_{\text{ure}}} \sigma_i^2 \leq -\sum_{i=1}^{N} y_i^2 + 2\sum_{i=1}^{N} \sigma_i^2.$$

Note also that

$$\|\theta\|^2 - \sum_{i=1}^{N_{\text{ure}}} \theta_i^2 + \sum_{i=1}^{N_{\text{ure}}} \sigma_i^2 = \|\theta\|^2 - \sum_{i=1}^{N_{\text{ure}}} y_i^2 + 2\sum_{i=1}^{N_{\text{ure}}} \sigma_i^2$$
$$+ 2\sum_{i=1}^{N_{\text{ure}}} \theta_i \sigma_i \xi_i + \sum_{i=1}^{N_{\text{ure}}} \sigma_i^2(\xi_i^2 - 1).$$

Therefore, combining this display and (4.33), we see that for any $N \geq 1$

$$(4.34) \quad \mathbf{E}_\theta R(\theta, N_{\text{ure}}) \leq R(\theta, N) + 2\mathbf{E}_\theta \sum_{i=1}^{N_{\text{ure}}} \theta_i \sigma_i \xi_i + \mathbf{E}_\theta \sum_{i=1}^{N_{\text{ure}}} \sigma_i^2(\xi_i^2 - 1).$$

In order to control the interference term $\mathbf{E}\sum_{i=1}^{N_{\text{ure}}} \theta_i \sigma_i \xi_i$, we use (4.31) with $\tilde{N} = N_{\text{ure}}$. This yields

$$\mathbf{E}_\theta \sum_{i=1}^{N_{\text{ure}}} \theta_i \sigma_i \xi_i \leq \alpha R(\theta, N) + \alpha \mathbf{E}_\theta \|\tilde{\theta}(N_{\text{ure}}) - \theta\|^2 - \alpha \mathbf{E}_\theta \sum_{i=1}^{N_{\text{ure}}} \sigma_i^2(\xi_i^2 - 1) + \frac{C\sigma_1^2}{\alpha^{4\beta+1}}.$$

Substituting this in (4.34), we have

$$(4.35) \quad \begin{aligned} \mathbf{E}_\theta R(\theta, N_{\text{ure}}) &\leq (1+2\alpha)R(\theta, N) + 2\alpha \mathbf{E}_\theta \|\tilde{\theta}(N_{\text{rhm}}) - \theta\|^2 \\ &\quad + (1-2\alpha)\mathbf{E}_\theta \sum_{i=1}^{N_{\text{ure}}} \sigma_i^2(\xi_i^2 - 1) + \frac{C\sigma_1^2}{\alpha^{4\beta+1}}. \end{aligned}$$

The last term in the above inequality can be controlled by Lemma 5, which gives that for any sufficiently small $\alpha > 0$ and $Q \in (1/2, (2\beta+1)/(4\beta+1)]$,

$$\mathbf{E}_\theta \sum_{i=1}^{N_{\text{ure}}} \sigma_i^2(\xi_i^2 - 1) \leq \alpha \sigma_1^2 \mathbf{E}_\theta \left(\sum_{i=1}^{N_{\text{ure}}} \frac{\sigma_i^4}{\sigma_1^4}\right)^Q + \frac{C\sigma_1^2}{\alpha^{1/(2Q-1)}}.$$

Let $Q = (2\beta + 1)/(4\beta + 1)$. Then by (2.6)

$$\left(\sum_{i=1}^{N_{\text{ure}}} \sigma_i^4\right)^Q \leq C\sigma_1^{4Q-2} \sum_{i=1}^{N_{\text{ure}}} \sigma_i^2,$$

and thus we obtain

$$(4.36) \qquad \mathbf{E}_\theta \sum_{i=1}^{N_{\text{ure}}} \sigma_i^2(\xi_i^2 - 1) \leq C\mu \mathbf{E}_\theta \sum_{i=1}^{N_{\text{ure}}} \sigma_i^2 + \frac{C\sigma_1^2}{\mu^{4\beta+1}}.$$

Finally, combining the above equation with (4.35), (4.36) and (4.32), we complete the proof. □

**Acknowledgments.** We would like to thank two anonymous referees, an Associate Editor and Co-Editor Jianqing Fan, whose constructive comments helped to improve the presentation of the paper.

CMI
39 RUE F. JOLIOT-CURIE
13453 MARSEILLE
CEDEX 13
FRANCE
E-MAIL: cavalier@cmi.univ-mrs.fr
golubev@cmi.univ-mrs.fr